\documentclass[11pt]{article}
\input epsf.tex

\usepackage{graphicx}
\usepackage{amsmath,amsthm,amsfonts,amscd,amssymb,comment,eucal,latexsym,mathrsfs}
\usepackage{stmaryrd}
\usepackage[all]{xy}

\usepackage{epsfig}

\usepackage[all]{xy}
\xyoption{poly}
\usepackage{fancyhdr}
\usepackage{wrapfig}
\usepackage{epsfig}

\theoremstyle{plain}
\newtheorem{thm}{Theorem}[section]

\newtheorem{conj}{Conjecture}

\theoremstyle{definition}

\theoremstyle{remark}
\newtheorem{remark}{Remark}
\newtheorem{example}{Example}

\advance \topmargin by -\headheight
\advance \topmargin by -\headsep
\evensidemargin \oddsidemargin
    \def\E{{\mathbb{E}}} \def\F{{\mathbb{F}}}        \def\N{{\mathbb{N}}}    \def\R{{\mathbb{R}}}  \def\T{{\mathbb{T}}}      \def\Z{{\mathbb{Z}}}

\def\cA{{\mathcal{A}}} \def\cB{{\mathcal{B}}}             \def\cO{{\mathcal{O}}} \def\cP{{\mathcal{P}}}     \def\cU{{\mathcal{U}}}     

        \def\tpsi{{\widetilde{\psi}}}        

                       \def\ux{{\underline{x}}} \def\uy{{\underline{y}}} 

\newcommand{\G}{\Gamma}

\newcommand{\Si}{\Sigma}
\newcommand{\eps}{\epsilon}

\renewcommand\d{\delta}

\renewcommand\k{\kappa}
\renewcommand\l{\lambda}

\newcommand\s{\sigma}

\newcommand\Aut{\operatorname{Aut}}

\newcommand\Hom{\operatorname{Hom}}

\newcommand\Prob{\operatorname{Prob}}

\newcommand\Sep{\operatorname{Sep}}

\newcommand\sym{\operatorname{sym}}

\def\cc{{\curvearrowright}}

\begin{document}
\title{A brief introduction to sofic entropy theory}
\author{Lewis Bowen\footnote{supported in part by NSF grant DMS-1500389 and a Simons Fellowship} \\ University of Texas at Austin}
\maketitle

\begin{abstract}
Sofic entropy theory is a generalization of the classical Kolmogorov-Sinai entropy theory to actions of a large class of non-amenable groups called sofic groups. This is a short introduction with a guide to the literature.
\end{abstract}


\noindent

\section{Introduction}

Classical entropy theory is concerned with single transformations of a topological or measure space. This can be generalized straightforwardly to actions of the lattice $\Z^d$. However, one encounters real difficulty in any attempt to generalize to actions of non-amenable groups. This short survey will begin with the free group of rank 2, $\F_2:=\langle a,b\rangle$. The Cayley graph $G=(V,E)$ of this group has vertex set $V=\F_2$ and directed edges $(g,ga), (g,gb)$ for $g\in \F_2$. It is a 4-regular tree. It is non-amenable because any finite subset $F \subset V$ has the property that if $\partial F$ is the set of edges $e \in E$ with one end in $F$ and one end outside of $F$ then $|\partial F| \ge 2|F|$. After understanding the special case of the free group (from a dynamicist's view), we will generalize to residually finite groups and sofic groups and then briefly survey recent developments; namely the classification of Bernoulli shifts, Bernoulli factors, Rokhlin entropy theory, algebraic dynamics and the geometry of model spaces.

We will not define amenability here (see \cite{MR3616077} for example). We will also not cover classical entropy theory. The interested reader is encouraged to consult one of the standard texts (e.g. \cite{petersen-book}). Other introductions and surveys on sofic entropy theory include \cite{MR3411529, MR3616077, gaboriau-sofic-survey, bowen-survey}.

\subsection{The Ornstein-Weiss factor}
In 1987, Ornstein and Weiss exhibited a curious example \cite{OW87}. To explain it, let $X:=(\Z/2\Z)^{\F_2}$ be the set of all maps $x:\F_2 \to \Z/2\Z$. This is a compact abelian group under pointwise addition. We can identify $X\times X$ with the group of maps from $\F_2 \to \Z/2\Z\times \Z/2\Z$. Define $\Phi:X \to X\times X$ by
$$\Phi(x)(g)=(x(g)-x(ga), x(g)-x(gb)).$$
This a surjective homomorphism.  It is also $\F_2$-equivariant where $\F_2$ acts on $X$ by
$$(fx)(g):=x(f^{-1}g).$$
And its kernel consists of the two constant maps. So it is a continuous, algebraic, 2-1 factor map from $(\Z/2\Z)^{\F_2}$ onto $(\Z/2\Z\times \Z/2\Z)^{\F_2}$. 

This appears to give a contradiction to entropy theory because if $K$ is a finite set then the entropy of the action of $\F_2$ on $K^{\F_2}$ should be $\log |K|$. So the Ornstein-Weiss factor map increases entropy! At the time of \cite{OW87} it was unknown whether or not the two actions $\F_2 \cc X$ and $\F_2 \cc X\times X$ could be measurably conjugate (with respect to Haar measure on $X$ and $X\times X$). We will show by Theorem \ref{thm:bernoulli1} below that they are not measurably conjugate. 

\section{Topological entropy for $\Z$-actions}

Here we will develop entropy theory for $\Z$-actions in a slightly non-traditional way which generalizes to actions of free groups. To begin, consider a homeomorphism $T:X \to X$ of a compact metric space $(X,d)$. A {\bf partial orbit} of length $n$ is a tuple of the form $(x,Tx,\ldots, T^{n-1}x) \in X^n$. Define a metric $d^{(n)}_\infty$ on $X^n$ by
$$d^{(n)}_\infty( \ux, \uy) = \max_i d(x_i,y_i)$$
where, for example $\ux=(x_1,\ldots, x_n)$. Then Rufus Bowen's definition of the topological entropy of $T$ is
$$h(T):=\sup_{\eps>0} \limsup_{n\to\infty} \frac{1}{n} \log \Sep_\eps( \{\textrm{length $n$ partial orbits}\}, d^{(n)}_\infty)$$
where if $S \subset X^n$ then $\Sep_\eps(S, d_\infty)$ denotes the maximum cardinality of an $\eps$-separated subset $Y \subset S$ \cite{MR0274707}. (Recall that a subset $Y \subset S$ is {\bf $\eps$-separated} if $d(y,z)>\eps$ for any $y,z\in Y$ with $y\ne z$).




Instead of counting partial orbits to compute entropy, we can count pseudo-orbits. To be precise, an $n$-tuple $\ux \in X^n$ is an {\bf $(n,\d)$-pseudo orbit} if 
$$\frac{1}{n} \sum_{i=1}^{n-1} d(Tx_i, x_{i+1}) < \d.$$
By Markov's inequality and continuity, any pseudo orbit contains a long subword that is close to a partial orbit. Hence
$$h(T)=\sup_{\eps>0}\inf_{\d>0} \limsup_{n\to\infty} \frac{1}{n} \log \Sep_\eps(\{ \textrm{$(n,\d)$ pseudo orbits}\}, d^{(n)}_\infty).$$
This is similar to Katok's treatment of entropy in \cite{MR573822} and follows from Kerr-Li's approach to entropy in \cite{kerr-li-variational}. 

It is of classical interest to count periodic orbits too. We will say that $\ux \in X^n$ is {\bf periodic with period $\le n$} if $Tx_i=x_{i+1}$ for all $1\le i \le n-1$ and $Tx_n=x_1$. The growth rate of periodic orbits is a lower bound for the entropy rate but in general they are not equal. To remedy this, let us consider pseudo-periodic orbits. To be precise, a $n$-tuple $\ux \in X^n$ is an {\bf $(n,\d)$-pseudo-periodic orbit} if
$$\frac{1}{n} \left( \sum_{i=1}^{n-1} d(Tx_i, x_{i+1}) + d(Tx_n, x_1)\right) < \d.$$
Since $(n,\d)$-pseudo-periodic orbits are $(n,\d)$-pseudo orbits and $(n,\d)$-pseudo orbits are $(n,\d + o_n(1))$-pseudo-periodic orbits, it follows that 
$$h(T)=\sup_{\eps>0}\inf_{\d>0} \limsup_{n\to\infty} \frac{1}{n} \log \Sep_\eps( \{\textrm{$(n,\d)$ pseudo-periodic orbits}\}, d^{(n)}_\infty).$$
Now it might look like we have not gained much by these observations since pseudo-periodic orbits are almost the same as pseudo orbits and the latter shadow partial orbits of only slightly less length. However, there is a conceptual advantage. This is because pseudo-periodic orbits can be thought of as maps from an external model, namely $\Z/n\Z$, to $X$ that approximate the dynamics. It is a very useful observation that we are not required to count only partial orbits or periodic orbits, both of which are too restrictive to generalize to actions of free groups.

\section{Topological sofic entropy for actions of free groups}

Suppose the free group $\F_2$ acts on a compact space $X$ by homeomorphisms. A {\bf periodic orbit} of this action consists of a finite set $V_0$, an action of $\F_2$ on $V_0$ and an $\F_2$-equivariant map $\phi:V_0 \to X$. 

It can be helpful to visualize the action of $\F_2$ on $V_0$ by making the {\bf action graph} $G_0=(V_0,E_0)$ whose edges consist of all pairs of the form $(v, a\cdot v)$ and $(v,b\cdot v)$ for $v\in V_0$. It is a directed graph in which every vertex has in-degree and out-degree 2. 

Given a finite subset $F \subset \F_2$ and $\d>0$, a {\bf $(V_0,\d,F)$-pseudo-periodic orbit} is a map $\phi:V_0 \to X$ that is approximately equivariant in the sense that
$$|V_0|^{-1} \sum_{v\in V_0} d( \phi(g\cdot v), g\cdot \phi(v)) < \d$$
for every $g\in F$. 

Now suppose that we fix a sequence $\Si:=\{\F_2 \cc V_i\}_{i=1}^\infty$ of actions of $\F_2$ on finite sets $V_i$. Tentatively, we will call the {\bf topological sofic entropy of $\F_2 \cc X$ with respect to $\Si$} the quantity
$$h_\Si(\F_2 \cc X) := \sup_{\eps>0}\inf_{\d>0} \inf_{F \Subset \F_2} \limsup_{n\to\infty} \frac{1}{n} \log \Sep_\eps( \{\textrm{$(V_n,\d,F)$ pseudo-periodic orbits}\}, d^{V_n}_\infty).$$
where $d^{V_n}_\infty$ is the metric on $X^{V_n}$ defined by
$$d^{V_n}_\infty(\phi,\psi) := \max_{v\in V_n} d(\phi(v),\psi(v)).$$

To state the next result, we need some terminology. Given a countable group $\G$ and continuous actions $\G \cc X$, $\G \cc Y$, an {\bf embedding} of $\G \cc X$ into $\G \cc Y$ is a continuous $\G$-equivariant injective map $\Phi:X\to Y$. If it is also surjective then it is a {\bf topological conjugacy}.

\begin{thm}
If $\F_2 \cc X$ embeds into $\F_2 \cc Y$ then for any $\Si$,
$$h_\Si(\F_2  \cc X) \le h_\Si(\F_2 \cc Y).$$
In particular, topological sofic entropy is a topological conjugacy invariant.
\end{thm}
The proof of this is straightforward, see \cite{kerr-li-variational, MR3616077} for details. (The definition of topological sofic entropy is due to Kerr-Li \cite{kerr-li-variational} which was inspired by my earlier work \cite{MR2552252}). 

\subsection{Examples}

\subsubsection{A boring example and asymptotic freeness}\label{boring}

Suppose that $V_n$ is a single point for all $n$. In this case, $h_\Si(\F_2 \cc X)$ is simply the logarithm of the number of fixed points of the action. While this is a topological conjugacy invariant, it is not what one usually means by entropy. To avoid this kind of example, we require that the actions $\{\F_2 \cc V_n\}_n$ are {\bf asymptotically free}. This means that for every nonidentity element $g\in \F_2$ 
\begin{eqnarray}\label{eqn:free}
\lim_{n\to\infty} |V_n|^{-1} \#\{v \in V_n:~ g\cdot v = v\} = 0.
\end{eqnarray}
A countable group $\G$ admits a sequence $\{\G \cc V_n\}_{n=1}^\infty$ of actions on finite sets satisfying asymptotic freeness if and only if $\G$ is {\bf residually finite}. We will come back to this point later. From now on, we assume the actions $\{\F_2 \cc V_n\}_n$ are asymptotically free.

\subsubsection{A curious example}\label{curious}

Let $X=\Z/2\Z$ and consider the action $\F_2 \cc X$ defined by $s \cdot x = x+1$ for $s\in \{a,b\}$. Let $\Si=\{\G \cc V_n\}_{n=1}^\infty$ be a sequence of actions on finite sets with the property that the corresponding action graphs $G_n:=(V_n,E_n)$ are bipartite. Let $V_n=P_n \sqcup Q_n$ be the bi-partition. Then define $\phi:V_n \to X$ by $\phi(P_n)=\{0\}$ and $\phi(Q_n)=\{1\}$. This map is a $(V_n,\d,F)$-pseudo-periodic orbit for all $\d,F$. So $h_\Si(\F_2 \cc X)\ge 0$. It can be shown that in fact any pseudo-periodic orbit must be close to either $\phi$ or $\phi+1$ which implies $h_\Si(\F_2 \cc X)=0$.

Next let $\Si'=\{\G \cc V'_n\}_{n=1}^\infty$ be a sequence of actions such that the corresponding action graphs $G'_n:=(V'_n,E'_n)$ are far from bi-partite. For example, it is known (and will be explained in \S \ref{f}) that if the action $\G \cc V'_n$ is chosen uniformly at random and $|V'_n| \to \infty$ as $n\to\infty$ then with high probability the action graphs will be far from bi-partite in the following sense. For small enough $\d>0$ and $F=\{a,b,a^{-1},b^{-1}\}$ there are no $(V'_n,\d,F)$-pseudo-periodic orbits. Since $\log(0)=-\infty$ this implies $h_{\Si'}(\G \cc X)=-\infty$.

This example shows (1) entropy depends on the choice of sequence $\Si$ and (2) it is possible for the entropy to be $-\infty$, even for very simple systems.

\subsubsection{The Ornstein-Weiss example revisited}

In \S \ref{sec:symbolic-top} below we will sketch a proof that the full shift action $\F_2 \cc K^{\F_2}$ has topological sofic entropy $\log |K|$ for any finite set $K$. So the Ornstein-Weiss factor map does indeed increase entropy. How can this happen? The answer is that if $\F_2 \cc X$ {\bf factors} onto $\F_2 \cc Y$ (meaning there is a continuous $\F_2$-equivariant surjective map $\Phi:X \to Y$) then, generally speaking, there is no way to ``lift'' pseudo-periodic orbits of the downstairs action $\F_2 \cc Y$ up to the source action $\F_2 \cc X$, even approximately. 

Let us see this in detail for the Ornstein-Weiss map. Suppose $\Si:=\{\F_2 \cc V_n\}_n$ is a sequence of actions on finite sets and form the action graphs $G_n=(V_n,E_n)$. Given any map $\psi:V_n \to \Z/2\Z \times \Z/2\Z$, we can define the {\bf pullback} or {\bf pullback name} of $\psi$ by 
$$\tpsi:V_n \to (\Z/2\Z\times \Z/2\Z)^{\F_2}, \quad \tpsi(v)(g)=\psi(g^{-1} \cdot v).$$
This map is a $(V_n,\d,F)$-pseudo-periodic orbit for the shift-action $\F_2 \cc (\Z/2\Z\times \Z/2\Z)^{\F_2}$ for any $\d>0$ and finite $F\subset \F_2$. In particular, this can be used to show that the entropy of $\F_2 \cc (\Z/2\Z \times \Z/2\Z)^{\F_2}$ is at least $\log 4$.

However most of these maps do not ``lift" via the Ornstein-Weiss map. To be precise, define
$$\Phi_n: (\Z/2\Z)^{V_n} \to (\Z/2\Z \times \Z/2\Z)^{V_n}, \quad \Phi_n(\psi)(v) = (\psi(v) - \psi(a^{-1}\cdot v), \psi(v) - \psi(b^{-1}\cdot v)).$$
This map is induced from the Ornstein-Weiss map. The point is that while the Ornstein-Weiss map is surjective, its finite approximations $\Phi_n$ are from surjective. This is obvious since the domain of $\Phi_n$ is exponentially smaller than  $(\Z/2\Z \times \Z/2\Z)^{V_n}$. Another argument is homological.

 Given $\psi:V_n \to \Z/2\Z \times \Z/2\Z$, let $\psi': E_n \to \Z/2\Z$ be the map defined by
$$\psi(v) = \left( \psi'(v,a^{-1}\cdot v), \psi'(v,b^{-1}\cdot v) \right).$$
Then $\psi$ is in the image of $\Phi_n$ if and only if $\psi'$ is a coboundary. So the ``reason'' the Ornstein-Weiss factor map increases entropy is that the $\Z/2\Z$-homology of the approximating graphs $G_n$ grows exponentially.  This observation generalizes: Gaboriau and Seward show in \cite{gaboriau-seward-2015} that if $\G$ is any sofic group and $k$ is a finite field then the sofic entropy of $\G \cc k^\G/k$ is at least $(1+\beta^1_{(2)}(\G))\log |k|$ where $\beta^1_{(2)}(\G)$ is the first $\ell^2$-Betti number of $\G$ and $k^\G/k$ is the quotient of $k^\G$ by the constant functions.

By contrast, any pseudo-periodic orbit of a $\Z$-action is close to a partial orbit of slightly less length. Partial orbits always lift. This explains why entropy is monotone decreasing for actions of $\Z$.

\section{Sofic groups}

You might have noticed that we have not used any special properties of free groups. In fact, the definition of topological sofic entropy stated above works for any residually finite group $\G$ in place of $\F_2$. Recall that $\G$ is {\bf residually finite} if there exists a decreasing sequence $\G \ge \G_1 \ge \G_2 \cdots$ such that each $\G_n$ is normal and finite-index in $\G$ and $\cap_n \G_n = \{e\}$ (this is equivalent to the previous definition in \S \ref{boring}). In this case, the sequence of actions $\G \cc \G/\G_n$ is asymptotically free and so the above definition of topological sofic entropy makes sense.

However, the actions $\G \cc V_n$ do not really have to be actions! To explain, let $\s_n:\G \to \sym(V_n)$ be a sequence of maps from $\G$ to the symmetric groups $\sym(V_n)$. We do not require these maps to be homomorphisms but we do require that they are {\bf asymptotically multiplicative} in the following sense: for every $g,h \in \G$ we require:
\begin{eqnarray}\label{eqn:homo}
\lim_{n\to\infty} |V_n|^{-1} \#\{v\in V_n:~\s_n(gh)v = \s_n(g)\s_n(h)v\}= 1.
\end{eqnarray}
We still require that they are asymptotically free, which means for every nonidentity $g\in \G$ 
\begin{eqnarray}\label{eqn:free2}
\lim_{n\to\infty} |V_n|^{-1} \#\{v \in V_n:~ \s_n(g)\cdot v = v\} = 0.
\end{eqnarray}
Any sequence $\Si=\{\s_n\}_{n=1}^\infty$ satisfying equations (\ref{eqn:homo}, \ref{eqn:free2}) is called a {\bf sofic approximation to $\G$} and $\G$ is called {\bf sofic} if it has a sofic approximation. The definition of topological sofic entropy given above makes sense for arbitrary sofic groups with respect to a sofic approximation $\Si$ once we replace $g \cdot v$ in the definition of a pseudo-periodic orbit with $\s_n(g)\cdot v$. 

Sofic groups were defined implicitly by Gromov in \cite{MR1694588}. They were given their name by Benjy Weiss in \cite{weiss-2000}. It is known that amenable groups and residually finite groups are sofic. The class of sofic groups is closed under a large number of group operations including passing to subgroups, direct limits, inverse limits, extensions by amenable groups, direct products, free products with amalgamation over amenable subgroups, graph products and wreathe products \cite{elek-szabo-2006, dykema-2014, paunescu-2011, elek-szabo-2011, CHR14, hayes-sale1}. By Malcev's Theorem, finitely generated linear groups are residually finite \cite{malcev-1940}. Since soficity is closed under direct limits, all countable linear groups are sofic.  Sofic groups solve special cases of a number of general conjectures including Connes Embedding Conjecture \cite{elek-szabo-2005}, the Determinant Conjecture \cite{elek-szabo-2005}, the Algebraic Eigenvalue Conjecture  \cite{MR2417890} and Gottschalk's Surjunctivity Conjecture (more on that later on). It is a major open problem whether all countable groups are sofic. Surveys on sofic groups include \cite{pestov-sofic-survey, pestov-kwiatkowska, capraro-lupini}.

\section{An application to Gottschalk's Surjunctivity Conjecture}

\begin{conj}\cite{MR0407821}
Let $k$ be a finite set, $\G$ a countable group and $\Phi:k^\G \to k^\G$ a continuous $\G$-equivariant map (where $k^\G$ is given the product topology). If $\Phi$ is injective then it is also surjective.
\end{conj}

This conjecture was proven to be true whenever $\G$ is sofic by Gromov \cite{MR1694588}. Another proof was given by Weiss \cite{weiss-2000} and then another by Kerr-Li \cite{kerr-li-variational}. Here is a sketch of Kerr-Li's proof: the sofic entropy of $\G \cc k^\G$ is $\log |k|$. However, the sofic entropy of any proper closed $\G$-invariant subset $X \subset k^\G$ is strictly less than $\log |k|$. Since entropy is a topological invariant, this proves the conjecture.

By the way, Gottschalk's conjecture implies Kaplansky's Direct Finiteness Conjecture which states: if $k$ is a finite field, $x,y$ are elements of the group ring $k\G$ and $xy=1$ then $yx=1$. To see the connection, observe that $x$ and $y$ induce linear $\G$-equivariant maps $\Phi_x, \Phi_y$ from $k^\G$ to itself. If $xy=1$ then $\Phi_x\Phi_y$ is the identity and so $\Phi_y$ is injective. If Gottschalk's conjecture holds then $\Phi_y$ must also be surjective and therefore $\Phi_x$ is its inverse. So $\Phi_y \Phi_x =1$, which implies $yx=1$. In fact this result holds for all fields $k$, even infinite fields because every field is embeddable into an ultraproduct of finite fields \cite{capraro-lupini}. When $\G$ is sofic, direct finiteness of $k\G$ also holds whenever $k$ is a division ring \cite{MR2089244} or a unital left Noetherian ring \cite{li-liang-mean-length}.

\section{Measure sofic entropy}

Before defining measure sofic entropy, let us revisit topological sofic entropy. The new notation will be useful in treating the measure case. Again, let $\G \cc X$ be an action by homeomorphisms on a compact metric space $(X,d)$ and fix a sofic approximation $\Si=\{\s_n\}_{n=1}^\infty$ where $\s_n: \G \to \sym(V_n)$. Let $\Omega( \s_n, \d,F) \subset X^{V_n}$ be the set of all $(\s_n,\d,F)$-pseudo-periodic orbits. To be precise, $\phi\in X^{V_n}$ is in $\Omega(\s_n,\d,F)$ if and only if
$$|V_n|^{-1} \sum_{v\in V_n} d( \phi(\s_n(g)\cdot v), g\cdot \phi(v)) < \d$$
for every $g\in F$. So $\Omega(\s_n,\d,F)$ depends implicitly on the action $\G \cc X$. The topological sofic entropy of $\G \cc X$ is defined by
$$h_\Si(\G \cc X) := \sup_{\eps>0}\inf_{\d>0} \inf_{F \Subset \G} \limsup_{n\to\infty} \frac{1}{n} \log \Sep_\eps( \Omega(\s_n, \d,F), d^{V_n}_\infty).$$

To define measure sofic entropy, we need a few more preliminaries. Let $\Prob(X)$ denote the space of Borel probability measures on $X$. It is a compact metrizable space with respect to the weak* topology which is defined by: a sequence of measures $\mu_n \in \Prob(X)$ converges to a measure $\mu$ if and only if: for every continuous function $f$ on $X$,
$$\int f~d\mu_n \to \int f~d\mu$$
as $n\to\infty$. 

Given a map $\phi:V_0 \to X$ (where $V_0$ is a finite set), the {\bf empirical distribution} of $\phi$ is the measure
$$P_\phi:= |V_0|^{-1} \sum_{v\in V_0} \d_v \in \Prob(X).$$

Given an open subset $\cO \subset \Prob(X)$, let $\Omega(\s_n,\d,F,\cO)$ be the set of all pseudo-periodic orbits $\phi \in \Omega(\s_n,\d,F)$ such that $P_\phi \in \cO$. Then the {\bf sofic entropy of a measure-preserving action $\G \cc (X,\mu)$ with respect to $\Si$} is
$$h_\Si(\G \cc (X,\mu)) := \sup_{\eps>0}\inf_{\d>0} \inf_{F \Subset \G} \inf_{\cO \ni \mu} \limsup_{n\to\infty} \frac{1}{n} \log \Sep_\eps( \Omega(\s_n, \d,F, \cO), d^{V_n}_\infty).$$

This is a measure-conjugacy invariant! For the proof see \cite{kerr-li-variational, MR3616077}. This definition is due to Kerr-Li \cite{kerr-li-variational} which was inspired by my earlier efforts \cite{MR2552252}. The proof is straightforward if the measure-conjugacy is a topological conjugacy. The general case is obtained by approximating a measure-conjugacy and its inverse by continuous maps. 

There is also a variational principle (due to Kerr-Li \cite{kerr-li-variational}):

\begin{thm}[Variational Principle]
$h(\G \cc X) = \sup_\mu h(\G \cc (X,\mu))$
where the sup is over all $\G$-invariant measures $\mu \in \Prob(X)$. If no such measures exist then $h(\G \cc X)=-\infty$.
\end{thm}

There is also a notion of sofic pressure and a corresponding variational principle \cite{chung-pressure}. See also \cite{zhang-local-variational} for local versions. Moreover, sofic entropy (both topological and measure) agrees with classical entropy whenever $\G$ is amenable \cite{kerr-li-sofic-amenable, MR2901354}. 

\section{Symbolic actions: the topological case}\label{sec:symbolic-top}

Let $(\cA,d_\cA)$ be a compact metric space. In most applications, $\cA$ is either finite or a torus (thought of as a compact abelian group). Let $\cA^\G$ be the space of all functions $x:\G \to \cA$ with the topology of pointwise convergence on finite sets. The group $\G$ acts on this space by $(g\cdot x)(f)  = x(g^{-1}f)$. 

Now suppose $X \subset \cA^\G$ is a closed $\G$-invariant subspace. There is a more convenient definition of the sofic entropy of $\G \cc X$ based on maps $\phi:V_n \to \cA$ (instead of maps $\phi:V_n \to X$). Given $\phi:V_n \to \cA$ and $v \in V_n$, define the {\bf pullback of $\phi$} by 
$$\Pi^{\s_n}_v(\phi) \in \cA^\G, \quad \Pi^{\s_n}_v(\phi)(g):=\phi( \s_n(g)^{-1}v).$$

Given an open neighborhood $\cU$ of $X$ in $\cA^\G$, let $\Omega'(\s_n, \d,\cU)$ be the set of all maps $\phi:V_n \to \cA$ such that 
$$|V_n|^{-1} \#\{v \in V_n:~\Pi^{\s_n}_v(\phi) \in \cU\} \ge 1-\d.$$
We call such a map a {\bf $(\s_n,\d,\cU)$-microstate} (this terminology is inspired by Voiculescu's free entropy \cite{MR1403924}).  We also call such a map a {\bf microstate} if the parameters are understood or intentionally left ambiguous. Then
$$h_\Si(\G \cc X) = \sup_{\eps>0}\inf_{\d>0} \inf_{\cU \supset X} \limsup_{n\to\infty} \frac{1}{n} \log \Sep_\eps( \Omega'(\s_n, \d,\cU), d^{V_n}_{\infty})$$
where (by abuse of notation) the metric $d^{V_n}_\infty$ on $\cA^{V_n}$ is defined by $d^{V_n}_\infty(\phi,\psi)=\max_{v\in V_n} d_\cA(\phi(v),\psi(v))$. 

The reason this works is that given any $\phi' \in \Omega'(\s_n,\d,\cU)$ we can find a pseudo-periodic orbit $\phi:V_n \to X$ such that $\phi(v)$ is close to $\Pi^{\s_n}_v(\phi')$ for most $v$. Conversely, given a pseudo-periodic orbit $\phi$ we can define $\phi':V_n \to \cA$ by $\phi'(v) = $ the projection of $\phi$ to the identity-coordinate. Then $\Pi^{\s_n}_v(\phi')$ will be close to $\phi(v)$ for most $v$. The measure case of this statement is proven in \cite{MR3542515}. The topological case follows from the variational principle.

If $\cA$ is finite then we can simplify further by setting $\eps=0$. To be precise,
$$h_\Si(\G \cc X) = \inf_{\d>0} \inf_{\cU \supset X} \limsup_{n\to\infty} \frac{1}{n} \log\# \Omega'(\s_n, \d,\cU).$$
The reason this works is that if $\phi,\psi \in \cA^{V_n}$ are any distinct elements then $d^{V_n}_\infty(\phi,\psi) \ge c>0$ where $c$ is the minimum distance between distinct elements of $\cA$.

Using this definition of sofic entropy, it is easy to check that $h_\Si(\G \cc \cA^\G) = \log |\cA|$. Indeed, this is true because $\Omega'(\s_n,\d,\cU)=\cA^{V_n}$. 

\subsection{Symbolic actions: the measure case}\label{sec:symbolic-measure}

Suppose $\mu \in \Prob(\cA^\G)$ is a probability measure preserved under the action. For any open neighborhood $\cO \ni \mu$, a {\bf $(\s_n,\cO)$-microstate} for $\mu$ is a map $\phi:V_n \to \cA$ such that its {\bf empirical measure}, defined by
$$P_\phi:= |V_n|^{-1} \sum_{v\in V_n} \d_{\Pi^{\s_n}_v(\phi)}$$
is contained in $\cO$. Let $\Omega'(\s_n,\cO)$ be the set of all $(\s_n,\cO)$-microstates for $\mu$. Then
$$h_\Si(\G \cc (\cA^\G,\mu))= \sup_{\eps>0}\inf_{\cO \ni \mu} \limsup_{n\to\infty} \frac{1}{n} \log \Sep_\eps( \Omega'(\s_n, \cO), d^{V_n}_{\infty}).$$
This approach is proven in \cite{MR3542515}. 

As in the topological case, we can simplify further if $\cA$ is finite by setting $\eps=0$ to obtain
$$h_\Si(\G \cc (\cA^\G,\mu))=\inf_{\cO \ni \mu} \limsup_{n\to\infty} \frac{1}{n} \log  \#\Omega'(\s_n, \cO).$$
The above is essentially the same as my original definition of sofic entropy in \cite{MR2552252}. 

\subsection{Bernoulli shifts}
Let $\k$ be a Borel probability measure on $\cA$. The {\bf Shannon entropy} of $\k$ is
$$H(\k):= - \sum_{a\in \cA} \k(\{a\}) \log \k(\{a\})$$
if $\k$ is supported on a countable set (and $0\log 0 :=0$). If $\k$ is not supported on a countable set then $H(\k):=+\infty$.

Let $\k^\G$ be the product measure on $\cA^\G$. The action $\G \cc (\cA^\G,\k^\G)$ is called the {\bf Bernoulli shift over $\G$ with base space $(\cA,\k$)}. 
\begin{thm}\label{thm:bernoulli1}
For any sofic approximation $\Si$,
$$h_\Si(\G\cc (\cA^\G,\k^\G)) = H(\k).$$
\end{thm} 
The finite entropy case of this was obtained in \cite{MR2552252} and the infinite entropy case is in \cite{MR2813530}. Here is a sketch in the special case in which $\cA$ is finite. 

\noindent {\bf The lower bound}. Let $\phi$ be a random map $V_n \to \cA$ with distribution $\k^{V_n}$. A second moment argument shows that, for any open neighborhood $\cO$ of $\k$, if $n$ is sufficiently large then with high probability $\phi$ is a $(\s_n,\cO)$-microstate for $\k^\G$. By the law of large numbers or the Shannon-McMillan Theorem, any subset $S \subset \cA^{V_n}$ with measure close to 1 has cardinality at least  $e^{|V_n| H(\k)-o(|V_n|)}$. This proves the lower bound. 

\noindent {\bf The upper bound}. Suppose that $\cO$ consists of all measures $\mu \in \Prob(\cA^\G)$ such that if $P:\cA^\G \to \cA$ denotes projection onto the identity coordinate then $\|P_*\mu - \k\|_{TV}<\eps$. Then the number of $(\s_n,\cO)$-microstates is approximately the multinomial
$$ |V_n|! \left(\prod_{a\in \cA}  \lfloor \k(\{a\}) V_n \rfloor ! \right)^{-1}$$
which, by Stirling's formula, is approximately $e^{|V_n| H(\k)+o(|V_n|)}$. This proves the upper bound.

\subsection{The $f$-invariant and RS-entropy}\label{f}
For this section, consider the special case in which $\G=\F_r=\langle s_1,\ldots, s_r\rangle$ is the rank $r$ free group. Instead of fixing a sofic approximation, set $V_n:=\{1,\ldots,n\}$ and let $\s_n:\F_r \to \sym(V_n)$ be a uniformly random homomorphism. The {\bf $f$-invariant or RS-entropy} of a measure-preserving action $\F_r \cc (X,\mu)$ is defined in the same way as sofic entropy except that ones takes an expected value before the logarithm:
$$f(\mu):=h^{RS}(\mu):=\sup_{\eps>0}\inf_{\d>0} \inf_{F \Subset \G} \inf_{\cO \ni \mu} \limsup_{n\to\infty} \frac{1}{n} \log \E_{\s_n}\big[ \Sep_\eps( \Omega(\s_n, \d,F, \cO), d^{V_n}_\infty)\big].$$
The RS stands for ``replica-symmetric'' to emphasize the analogy with the corresponding notions in the literature  in statistical physics and theoretical computer science \cite{MR2643563}. In the special case in which $\mu$ is a shift invariant measure on $\cA^{\F_r}$ and $\cA$ is finite, the definition reduces to
$$f(\mu)=h^{RS}(\mu)=\inf_{\cO \ni \mu} \limsup_{n\to\infty} \frac{1}{n} \log \E_{\s_n}[ \#\Omega'(\s_n, \cO)].$$
Alternatively, let $\cP$ be a countable measurable partition of $X$. Its Shannon entropy is defined by
\begin{eqnarray}\label{eqn:Shannon}
H_\mu(\cP):= -\sum_{P\in \cP} \mu(P)\log\mu(P).
\end{eqnarray}
Define
$$F_\mu(\cP):= -(2r-1)H_\mu(\cP) + \sum_{i=1}^r H_\mu(\cP \vee s_i\cP)$$
where $\cP \vee s_i\cP$ is the smallest partition refining both $\cP$ and $s_i\cP$. For any finite $W \subset \F_r$, let $\cP^W$ be the smallest partition containing $w\cP$ for all $w\in W$.
\begin{thm}\cite{bowen-entropy-2010b}
$f(\mu)=\inf_{R>0} F_\mu(\cP^{B(R)})$ where $B(R) \subset \F_r$ is the ball of radius $R>0$ with respect to the word metric and $\cP$ is any generating partition for the action such that $H_\mu(\cP)<\infty$. (A partition $\cP$ is {\bf generating} if the smallest complete $\G$-invariant $\s$-sub-algebra containing $\cP$ consists of all measurable sets).
\end{thm}
The theorem above was taken as the definition of $f$ in \cite{bowen-annals-2010} where the $f$-invariant was first proven to be a measure-conjugacy invariant without using model spaces. Conditional on the existence of a finite generating partition, the $f$-invariant is additive under direct products, it satisfies an ergodic decomposition formula \cite{seward-ergodic}, a subgroup formula \cite{seward-subgroup}, an Abraham-Rokhlin formula \cite{MR2736889} and a (restricted) Yuzvinskii addition formula \cite{MR3163025}. Sofic entropy does not, in general, satisfy such formulas \cite{bowen-survey}. 

\begin{example}
The $f$-invariant of any action $\G$ on a finite set $X$ is $-(r-1)\log|X|$. So $\E_{\s_n}[ \#\Omega'(\s_n, \cO)] \approx |X|^{-(r-1)|V_n|}$ (for $\cO$ small and $n$ large). Since $ \#\Omega'(\s_n, \cO) \in \{0,1,2,\ldots\}$, with high probability $\Omega'(\s_n,\cO)$ is empty (if $\cO$ is small enough, $r\ge 2$ and $|X|\ge 2$). This explains why, as claimed in \S \ref{curious}, if $\G \cc V'_n$ is uniformly random then with high probability the graphs $G'_n=(V'_n,E'_n)$ are far from bipartite in the sense that there are no pseudo-periodic orbits (or microstates) for the action $\G \cc \Z/2\Z$ where each generator acts nontrivially. 
\end{example}

The paper \cite{MR2736889} defines Markov processes over the free group and shows that, for such processes, $F_\mu(\cP)=f(\mu)$. In particular, it is very easy to compute. It can be shown that some mixing Markov processes (e.g., the Ising process with small transition probability and free boundary conditions) have negative $f$-invariant. These cannot be isomorphic to Bernoulli shifts. By contrast, all mixing Markov processes over the integers are isomorphic to Bernoulli shifts. It is a major open problem to classify (mixing) Markov processes over a free group up to measure-conjugacy.

\section{Classification of Bernoulli shifts}

Theorem \ref{thm:bernoulli1} shows that if $\G$ is sofic then Bernoulli shifts with different base space entropies are not measurably conjugate. Surprisingly, the converse is true even without soficity:
\begin{thm}
If $\G$ is any countably infinite group and $(\cA,\k)$, $(\cB,\l)$ are two probability spaces with equal Shannon entropies $H(\k)=H(\l)$ then $\G \cc (\cA,\k)^\G$ is measurably conjugate to $\G \cc (\cB,\l)^\G$.
\end{thm}

\begin{remark}
The special case in which $\G=\Z$ is Ornstein's famous theorem of 1970 \cite{ornstein-1970a}. Stepin observed that if $\G$ is any countable group containing a copy of $\Z$ then the theorem above holds for $\G$ because one can build an isomorphism for $\G$-actions from an isomorphism for $\Z$-actions, coset-by-coset \cite{stepin-1975}. Ornstein-Weiss extended Ornstein's Theorem to all countably infinite amenable groups through the technology of quasi-tilings \cite{OW80}. I showed in \cite{MR2931910} that the above theorem is true whenever the supports of $\k$ and $\l$ each contain more than 2 elements. The proof is essentially a ``measurable version'' of Stepin's trick. The last remaining case (when say $|\cA|=2$) has been handled recently in soon-to-be-published work of Brandon Seward.

\end{remark}

\section{Bernoulli factors}

\begin{thm}\label{thm:factors}
Let $\G$ be any non-amenable group. Then every nontrivial Bernoulli shift over $\G$ factors onto every other nontrivial Bernoulli shift.
\end{thm}

This theorem is obtained in \cite{bowen-ri}. The special case of the free group $\F_2$ was handled in \cite{bowen-ornstein-2011} using the Ornstein-Weiss map and Sinai's Factor Theorem (for actions of $\Z$). It immediately follows for any group containing a copy of $\F_2$ since we can build the factor map coset-by-coset. In \cite{ball-factors1}, it was shown that if $\G$ is any non-amenable group then there is some Bernoulli shift over a finite base space that factors onto all other Bernoulli shifts. The argument used a rudimentary form of the Gaboriau-Lyons Theorem (before that theorem existed)  which states that if $\G$ is any non-amenable group then there is some Bernoulli shift (over a finite base space) $\G \cc (\cA,\k)^\G$ and an ergodic essentially free action of the free group $\F_2 \cc \cA^\G$ such that the orbits of the free group action are contained in the $\G$-orbits \cite{gaboriau-lyons}. We can then view this free group action as being essentially like having a free subgroup of $\G$ and build the factor map coset-by-coset as before. The main new result of \cite{bowen-ri} is that the Gaboriau-Lyons Theorem holds for arbitrary Bernoulli shifts. Theorem \ref{thm:factors} then follows from an argument similar to \cite{ball-factors1}.

\section{Rokhlin entropy}

A measurable partition $\cP$ of a measure space $(X,\mu)$ is {\bf generating} for an action $\G \cc (X,\mu)$ if the smallest $\G$-invariant sigma-algebra containing $\cP$ consists of all measurable sets (modulo sets of measure zero). 
The Rokhlin entropy of an ergodic action $\G \cc (X,\mu)$ is the infimum of the Shannon entropies of generating partitions (the Shannon entropy is defined by (\ref{eqn:Shannon})). In the special case that $\G=\Z$, Rokhlin proved that this agrees with Kolmogorov-Sinai entropy \cite{rohlin-lectures}. A modern proof of this, that holds for all amenable $\G$, is in \cite{seward-tucker-drob}. 

Rokhlin entropy is an upper bound to sofic entropy. It is unknown whether they are equal, conditioned on the sofic entropy not being minus infinity. For example, this is unknown even for principal algebraic actions (see \S \ref{principal}). Unfortunately, sofic entropy is the only known lower bound for Rokhlin entropy; hence we do not even know how to compute the Rokhlin entropy of Bernoulli shifts, except when the group is sofic, in which case the Rokhlin entropy equals the Shannon entropy of the base. Indeed, it is shown in \cite{seward-kreiger-2} that if the Rokhlin entropy of Bernoulli shifts is positive (for all groups) then Gottschalk's conjecture holds for all groups. It is also known that Rokhlin entropy equals sofic entropy for Gibbs measures satisfying a strong spatial mixing condition \cite{alpeev-percolation, austin-percolation}.

Rokhlin entropy has mainly been used as a hypothesis rather than a conclusion. There are two main theorems of this form; they  generalize Krieger's Generator Theorem and Sinai's Factor Theorem:

\begin{thm}\cite{seward-kreiger-1}
If $\G \cc (X,\mu)$ is ergodic and has Rokhlin entropy $< \log(n)$ for some integer $n>1$ then there exists a generating partition for the action with $n$ parts.
\end{thm}
This Theorem is the simplest version of a large variety of far more refined results contained in \cite{seward-kreiger-1, seward-kreiger-2, alpeev-seward}. 

\begin{thm}\cite{seward-sinai-30}
If $\G \cc (X,\mu)$ is ergodic and has positive Rokhlin entropy then it factors onto a Bernoulli shift.
\end{thm}

\section{Algebraic actions}

An {\bf algebraic action} is an action of a countable group $\G$ on a compact group $X$ by group-automorphisms. The main problem is to relate algebraic or analytic properties of the image of $\G$ in $\Aut(X)$ with purely dynamical properties. The study of single automorphisms (that is, $\G=\Z$) goes back at least to Yuzvinskii \cite{MR0201560, MR0194588} and R. Bowen \cite{MR0274707}. The special case $\G=\Z^d$ was studied intensively in the 80's and 90's \cite{schmidt-book}. Here we will highlight a few recent achievements and open problems  extending classical results to the realm of sofic group actions.

\subsubsection{Topological versus measure entropy}
Yuzvinskii and R. Bowen showed that, when $\G=\Z$, the topological and measure-entropy of an algebraic action $\G \cc X$ agree (where the measure on $X$ is Haar measure) \cite{MR0194588,MR0274707}. This was extended to amenable $\G$ by Deninger \cite{deninger-2006}. It was an open problem since 2011 whether this result could be extended to sofic groups. There were computations of entropy showing that it was true in a number of special cases \cite{kerr-li-variational, MR2794944, MR2956925, hayes-fk-determinants} but these all proceeded by computing the topological entropy and the measure entropy separately in terms of analytic data and then showing their equality. So it is astonishing that just recently Ben Hayes proved under a mild hypothesis on the actions that the topological entropy agrees with the measure entropy \cite{hayes-doubly-quenched}. The proof uses Austin's lde-sofic-entropy \cite{MR3542515}. 

\subsubsection{Principal algebraic actions}\label{principal}
 Let $f$ be an element of the integer group ring $\Z\G$ and consider the principal ideal $\Z\G f \subset \Z\G$ generated by $f$. Then $\Z\G/\Z\G f$ is a countable abelian group that $\G$ acts on by automorphisms (namely the action is $g (x+\Z\G f) := gx + \Z\G f$). Let $X_f:=\Hom(\Z\G/\Z\G f, \R/\Z)$ be the Pontryagin dual. This is a compact abelian group under pointwise addition. Moreover, the action of $\G$ on $\Z\G/\Z\G f$ induces an action of $\G$ on $X_f$ by automorphisms.

In the special case in which $\Z=\G$, we can write $f$ as $f=\sum_{i=-m}^n c_i x^i$ for some coefficients $c_i \in \Z$ and some $0\le m,n$. After multiplying $f$ by some $x^k$ we may assume that $m=0$ and $c_0\ne 0$. This change does not affect the dynamics.   So we can think of $f$ as a polynomial. Yuzvinskii and R. Bowen showed that the entropy of $\Z$ acting on $X_f$ is the log-Mahler measure of $f$. This result was extended to $\Z^d$ by Lind-Schmidt-Ward \cite{LSW-1990}. Chris Deninger observed that the Fuglede-Kadison determinant generalizes Mahler measure to non-abelian $\G$ and conjectured that, when $\G$ is amenable,   $h(\G \cc X_f)$ is the log of the Fuglede-Kadison determinant of $f$ \cite{deninger-2006}. Special cases were handled in \cite{deninger-2006, deninger-schmidt} before the general case was completed by Hanfeng Li in \cite{MR2925385}.

The first result in the setting of a non-amenable acting group $\G$ was the case of expansive principal algebraic actions of residually finite groups \cite{MR2794944}. This was extended in \cite{MR2956925} to some non-expansive actions (harmonic mod 1 points). Then in a stunning breakthrough Hayes proved that for an arbitrary sofic $\G$ and arbitrary $f \in \Z\G$ either $f$ is not injective as a convolution operator on $\ell^2(\G)$ (in which case the sofic entropy is infinite) or it is injective and the sofic entropy equals the log of the Fuglede-Kadison determinant of $f$ \cite{hayes-fk-determinants}. 

\subsubsection{Yuzvinskii's addition formula}

 Suppose that $N \vartriangleleft X$ is a closed $\G$-invariant normal subgroup. We say the {\bf addition formula} holds for $(\G \cc X,N)$ if the entropy of $\G \cc X$ equals the sum of the entropy of $\G \cc N$ with the entropy of $\G \cc X/N$. In the special case in which $\G=\Z$, this result is due to Yuzvinskii  \cite{MR0201560, MR0194588}. It was extended to $\Z^d$ in  \cite{LSW-1990} and to arbitrary amenable groups in \cite{MR2925385} (and independently in unpublished work of Lind-Schmidt). It is an important structural result which when combined with the principal algebraic case yields a general procedure for computing entropy of algebraic actions satisfying mild hypotheses \cite{schmidt-book}. Using it, Li and Thom relate entropy to $L^2$-torsion \cite{MR3110799} thereby obtaining new results about $L^2$-torsion and algebraic dynamics.
 
 The Ornstein-Weiss example shows that addition formulas fail for sofic entropy in general. Indeed, it has been shown in \cite{bartholdi-kielak} that for any non-amenable $\G$ and any field $k$, there exists an embedding $(k\G)^n \to (k\G)^{n-1}$ for some $n\in \N$ (as $k\G$-modules). Taking the dual and setting $k$ equal to a finite field gives a contradiction to the addition formula. 
 
However, the $f$-invariant satisfies the addition formula when $X$ is totally disconnected and satisfies some technical hypothesis \cite{MR3163025}. The reason this does not contradict the previous paragraph is that the $f$-invariant, unlike sofic entropy, can take finite negative values. It is an open problem whether the addition formula holds for the $f$-invariant in general.

\subsubsection{Pinsker algebra}
In recent work,  Hayes shows that under mild hypotheses, the outer Pinsker algebra of an algebraic action of a sofic group is algebraic; that is, it comes from an invariant closed normal subgroup \cite{hayes-relative-entropy}. To explain a little more, the outer sofic entropy of a factor is the growth rate of the number of microstates that lift to the source action. The outer Pinsker algebra is the maximal $\s$-sub-algebra such that the corresponding factor has zero outer sofic entropy. It follows that, in order, to prove an algebraic action has completely positive entropy (CPE), it is sufficiently to check all of its algebraic factors.

In \cite{schmidt-book} it is shown that if $\G=\Z^d$ then all CPE algebraic actions are Bernoulli. Could this be true more generally? Even for amenable groups, this problem is open.

\section{Geometry of model spaces}

Recall from \S \ref{sec:symbolic-measure} that
$$h_\Si(\G \cc (\cA^\G,\mu))= \sup_{\eps>0}\inf_{\cO \ni \mu} \limsup_{n\to\infty} \frac{1}{n} \log \Sep_\eps( \Omega'(\s_n, \cO), d^{V_n}_{\infty}).$$
The spaces $\Omega'(\s_n,\cO)$ are called {\bf model spaces}. We consider them with the metrics $d_1^{V_n}$ defined by
$$d_1^{V_n}(\phi, \psi):= \frac{1}{|V_n|} \sum_{v\in V_n} d(\phi(v),\psi(v)).$$
In the special case in which $\cA$ is finite and $d(x,y) \in \{0,1\}$ for all $x,y \in \cA$, $d_1^{V_n}$ is the normalized Hamming metric. The asymptotic {\it geometry} of these model spaces can be used to define new invariants. 

In \cite{MR3543677} Tim Austin introduced a notion of asymptotic coarse connectedness for model spaces that depends on a choice of sofic approximation $\Si$. He calls this notion {\bf connected model spaces rel $\Si$}. For fixed $\Si$, it is a measure-conjugacy invariant.  He shows that Bernoulli shifts have connected model spaces rel $\Si$ (for any $\Si$). On the other hand, if $\G$ is residually finite, has property (T) and $\T$ is the 1-torus then there is a sofic approximation $\Si$ such that the model spaces for the action $\G \cc \T^\G/\T$ are not asymptotically coarsely connected. In particular, this action is a factor of a Bernoulli shift that is not isomorphic to a Bernoulli shift (by contrast, all factors of Bernoulli shifts of $\Z$-actions are Bernoulli \cite{ornstein-1970c}). This was known earlier from work of Popa-Sasyk \cite{popa-sasyk} and the example is the same, but the proof is different since it goes through this new measure-conjugacy invariant. 

In forthcoming work, I will generalize model-connectedness to asymptotic coarse homological invariants. These new invariants will be applied to show that there are Markov processes over the free group that do not have the Weak Pinsker Property. By contrast, Austin recently proved that all processes over $\Z$ have the Weak Pinsker Property \cite{austin-wp}.

In \cite{MR3542515} Austin defines a notion of convergence for sequences of probability measures $\mu_n$ on $\cA^{V_n}$ (with respect to $\Si$) and uses this to define new notions of sofic entropy. One of these notions was called ``doubly-quenched sofic entropy'' in \cite{MR3542515} but has now been renamed to ``locally doubly empirical entropy'' to avoid conflict with physics notions. This particular version of sofic entropy is additive under direct products. Under very mild hypotheses, it agrees with the power-stabilized entropy which is defined by
$$h^{PS}_\Si(\G \cc (X,\mu)) = \lim_{n\to\infty} \frac{1}{n} h^{PS}_\Si(\G \cc (X,\mu)^n)$$
where $\G \cc X^n$ acts diagonally $g(x_1,\ldots, x_n):=(gx_1,\ldots, gx_n)$. Thus, it seems to be the `right' version of sofic entropy if one requires additivity under direct products. 
 
{\bf Acknowledgements}. Thanks to Tim Austin and Andrei Alpeev for helpful comments.


\bibliography{biblio}

\def\cprime{$'$} \def\cprime{$'$} \def\cprime{$'$}
  \def\cfudot#1{\ifmmode\setbox7\hbox{$\accent"5E#1$}\else
  \setbox7\hbox{\accent"5E#1}\penalty 10000\relax\fi\raise 1\ht7
  \hbox{\raise.1ex\hbox to 1\wd7{\hss.\hss}}\penalty 10000 \hskip-1\wd7\penalty
  10000\box7} \def\cprime{$'$} \def\cprime{$'$} \def\cprime{$'$}
  \def\cprime{$'$} \def\cprime{$'$} \def\cprime{$'$} \def\cprime{$'$}
\begin{thebibliography}{Bow17b}

\bibitem[Alp17]{alpeev-percolation}
Andrei Alpeev.
\newblock Random ordering formula for sofic and rokhlin entropy of gibbs
  measures.
\newblock {\em preprint}, 2017.

\bibitem[AP17]{austin-percolation}
Tim Austin and Moumanti Podder.
\newblock Gibbs measures over locally tree-like graphs and percolative entropy
  over infinite regular trees.
\newblock {\em preprint}, 2017.

\bibitem[AS16]{alpeev-seward}
Andrei Alpeev and Brandon Seward.
\newblock Krieger's finite generator theorem for ergodic actions of countable
  groups {III}.
\newblock {\em preprint}, 2016.

\bibitem[Aus16a]{MR3542515}
Tim Austin.
\newblock Additivity properties of sofic entropy and measures on model spaces.
\newblock {\em Forum Math. Sigma}, 4:e25, 79, 2016.

\bibitem[Aus16b]{MR3543677}
Tim Austin.
\newblock The {G}eometry of {M}odel {S}paces for {P}robability-{P}reserving
  {A}ctions of {S}ofic {G}roups.
\newblock {\em Anal. Geom. Metr. Spaces}, 4:Art. 6, 2016.

\bibitem[Aus17]{austin-wp}
Tim Austin.
\newblock Measure concentration and the weak pinsker property.
\newblock {\em preprint}, 2017.

\bibitem[Bal05]{ball-factors1}
Karen Ball.
\newblock Factors of independent and identically distributed processes with
  non-amenable group actions.
\newblock {\em Ergodic Theory Dynam. Systems}, 25(3):711--730, 2005.

\bibitem[BG14]{MR3163025}
Lewis Bowen and Yonatan Gutman.
\newblock A {J}uzvinskii addition theorem for finitely generated free group
  actions.
\newblock {\em Ergodic Theory Dynam. Systems}, 34(1):95--109, 2014.

\bibitem[BK17]{bartholdi-kielak}
L.~Bartholdi and D.~Kielak.
\newblock Amenability of groups is characterized by {M}yhill?s theorem.
\newblock {\em J. Eur. Math. Soc. to appear. arXiv:1605.09133v2.}, 2017.

\bibitem[BL12]{MR2956925}
Lewis Bowen and Hanfeng Li.
\newblock Harmonic models and spanning forests of residually finite groups.
\newblock {\em J. Funct. Anal.}, 263(7):1769--1808, 2012.

\bibitem[Bow71]{MR0274707}
Rufus Bowen.
\newblock Entropy for group endomorphisms and homogeneous spaces.
\newblock {\em Trans. Amer. Math. Soc.}, 153:401--414, 1971.

\bibitem[Bow10a]{bowen-entropy-2010b}
Lewis Bowen.
\newblock The ergodic theory of free group actions: entropy and the
  {$f$}-invariant.
\newblock {\em Groups Geom. Dyn.}, 4(3):419--432, 2010.

\bibitem[Bow10b]{MR2552252}
Lewis Bowen.
\newblock Measure conjugacy invariants for actions of countable sofic groups.
\newblock {\em J. Amer. Math. Soc.}, 23(1):217--245, 2010.

\bibitem[Bow10c]{MR2736889}
Lewis Bowen.
\newblock Non-abelian free group actions: {M}arkov processes, the
  {A}bramov-{R}ohlin formula and {Y}uzvinskii's formula.
\newblock {\em Ergodic Theory Dynam. Systems}, 30(6):1629--1663, 2010.

\bibitem[Bow10d]{bowen-annals-2010}
Lewis~Phylip Bowen.
\newblock A measure-conjugacy invariant for free group actions.
\newblock {\em Ann. of Math. (2)}, 171(2):1387--1400, 2010.

\bibitem[Bow11a]{MR2794944}
Lewis Bowen.
\newblock Entropy for expansive algebraic actions of residually finite groups.
\newblock {\em Ergodic Theory Dynam. Systems}, 31(3):703--718, 2011.

\bibitem[Bow11b]{bowen-ornstein-2011}
Lewis Bowen.
\newblock Weak isomorphisms between {B}ernoulli shifts.
\newblock {\em Israel J. Math.}, 183:93--102, 2011.

\bibitem[Bow12a]{MR2931910}
Lewis Bowen.
\newblock Every countably infinite group is almost {O}rnstein.
\newblock In {\em Dynamical systems and group actions}, volume 567 of {\em
  Contemp. Math.}, pages 67--78. Amer. Math. Soc., Providence, RI, 2012.

\bibitem[Bow12b]{MR2901354}
Lewis Bowen.
\newblock Sofic entropy and amenable groups.
\newblock {\em Ergodic Theory Dynam. Systems}, 32(2):427--466, 2012.

\bibitem[Bow17a]{bowen-survey}
Lewis Bowen.
\newblock Examples in the entropy theory of countable group actions.
\newblock {\em submitted}, 2017.

\bibitem[Bow17b]{bowen-ri}
Lewis Bowen.
\newblock Finitary random interlacements and the gaboriau-lyons problem.
\newblock {\em submitted}, 2017.

\bibitem[CHR14]{CHR14}
Laura Ciobanu, Derek~F. Holt, and Sarah Rees.
\newblock Sofic groups: graph products and graphs of groups.
\newblock {\em Pacific J. Math.}, 271(1):53--64, 2014.

\bibitem[Chu13]{chung-pressure}
Nhan-Phu Chung.
\newblock Topological pressure and the variational principle for actions of
  sofic groups.
\newblock {\em Ergodic Theory Dynam. Systems}, 33(5):1363--1390, 2013.

\bibitem[CL15]{capraro-lupini}
Valerio Capraro and Martino Lupini.
\newblock {\em Introduction to {S}ofic and hyperlinear groups and {C}onnes'
  embedding conjecture}, volume 2136 of {\em Lecture Notes in Mathematics}.
\newblock Springer, Cham, 2015.
\newblock With an appendix by Vladimir Pestov.

\bibitem[Den06]{deninger-2006}
Christopher Deninger.
\newblock Fuglede-{K}adison determinants and entropy for actions of discrete
  amenable groups.
\newblock {\em J. Amer. Math. Soc.}, 19(3):737--758 (electronic), 2006.

\bibitem[DKP14]{dykema-2014}
Ken Dykema, David Kerr, and Mika{\"e}l Pichot.
\newblock Sofic dimension for discrete measured groupoids.
\newblock {\em Trans. Amer. Math. Soc.}, 366(2):707--748, 2014.

\bibitem[DM10]{MR2643563}
Amir Dembo and Andrea Montanari.
\newblock Gibbs measures and phase transitions on sparse random graphs.
\newblock {\em Braz. J. Probab. Stat.}, 24(2):137--211, 2010.

\bibitem[DS07]{deninger-schmidt}
Christopher Deninger and Klaus Schmidt.
\newblock Expansive algebraic actions of discrete residually finite amenable
  groups and their entropy.
\newblock {\em Ergodic Theory Dynam. Systems}, 27(3):769--786, 2007.

\bibitem[ES04]{MR2089244}
G{\'a}bor Elek and Endre Szab{\'o}.
\newblock Sofic groups and direct finiteness.
\newblock {\em J. Algebra}, 280(2):426--434, 2004.

\bibitem[ES05]{elek-szabo-2005}
G{\'a}bor Elek and Endre Szab{\'o}.
\newblock Hyperlinearity, essentially free actions and {$L^2$}-invariants.
  {T}he sofic property.
\newblock {\em Math. Ann.}, 332(2):421--441, 2005.

\bibitem[ES06]{elek-szabo-2006}
G{\'a}bor Elek and Endre Szab{\'o}.
\newblock On sofic groups.
\newblock {\em J. Group Theory}, 9(2):161--171, 2006.

\bibitem[ES11]{elek-szabo-2011}
G{\'a}bor Elek and Endre Szab{\'o}.
\newblock Sofic representations of amenable groups.
\newblock {\em Proc. Amer. Math. Soc.}, 139(12):4285--4291, 2011.

\bibitem[Gab16]{gaboriau-sofic-survey}
Damien Gaboriau.
\newblock Sofic entropy, after lewis bowen, david kerr and hanfeng li.
\newblock {\em arXiv preprint arXiv:1607.06454}, 2016.

\bibitem[GL09]{gaboriau-lyons}
Damien Gaboriau and Russell Lyons.
\newblock A measurable-group-theoretic solution to von {N}eumann's problem.
\newblock {\em Invent. Math.}, 177(3):533--540, 2009.

\bibitem[Got73]{MR0407821}
Walter Gottschalk.
\newblock Some general dynamical notions.
\newblock pages 120--125. Lecture Notes in Math., Vol. 318, 1973.

\bibitem[Gro99]{MR1694588}
M.~Gromov.
\newblock Endomorphisms of symbolic algebraic varieties.
\newblock {\em J. Eur. Math. Soc. (JEMS)}, 1(2):109--197, 1999.

\bibitem[GS15]{gaboriau-seward-2015}
Damien Gaboriau and Brandon Seward.
\newblock Cost, $\ell^2$-betti numbers, and the sofic entropy of some algebraic
  actions.
\newblock {\em preprint}, 2015.

\bibitem[Hay16a]{hayes-doubly-quenched}
Ben Hayes.
\newblock Doubly quenched convergence and the entropy of algebriac actions of
  sofic groups.
\newblock {\em submitted}, 2016.

\bibitem[Hay16b]{hayes-fk-determinants}
Ben Hayes.
\newblock Fuglede-{K}adison determinants and sofic entropy.
\newblock {\em Geom. Funct. Anal.}, 26(2):520--606, 2016.

\bibitem[Hay16c]{hayes-relative-entropy}
Ben Hayes.
\newblock Relative entropy and the {P}insker product formula for sofic groups.
\newblock {\em submitted}, 2016.

\bibitem[HS16]{hayes-sale1}
Ben Hayes and Andrew Sale.
\newblock The wreath product of two sofic groups is sofic.
\newblock {\em submitted}, 2016.

\bibitem[Juz65a]{MR0201560}
S.~A. Juzvinski{\u\i}.
\newblock Metric properties of automorphisms of locally compact commutative
  groups.
\newblock {\em Sibirsk. Mat. \v Z.}, 6:244--247, 1965.

\bibitem[Juz65b]{MR0194588}
S.~A. Juzvinski{\u\i}.
\newblock Metric properties of the endomorphisms of compact groups.
\newblock {\em Izv. Akad. Nauk SSSR Ser. Mat.}, 29:1295--1328, 1965.

\bibitem[Kat80]{MR573822}
A.~Katok.
\newblock Lyapunov exponents, entropy and periodic orbits for diffeomorphisms.
\newblock {\em Inst. Hautes \'Etudes Sci. Publ. Math.}, (51):137--173, 1980.

\bibitem[KL11a]{MR2813530}
David Kerr and Hanfeng Li.
\newblock Bernoulli actions and infinite entropy.
\newblock {\em Groups Geom. Dyn.}, 5(3):663--672, 2011.

\bibitem[KL11b]{kerr-li-variational}
David Kerr and Hanfeng Li.
\newblock Entropy and the variational principle for actions of sofic groups.
\newblock {\em Invent. Math.}, 186(3):501--558, 2011.

\bibitem[KL13]{kerr-li-sofic-amenable}
David Kerr and Hanfeng Li.
\newblock Soficity, amenability, and dynamical entropy.
\newblock {\em Amer. J. Math.}, 135(3):721--761, 2013.

\bibitem[KL16]{MR3616077}
David Kerr and Hanfeng Li.
\newblock {\em Ergodic theory}.
\newblock Springer Monographs in Mathematics. Springer, Cham, 2016.
\newblock Independence and dichotomies.

\bibitem[Li12]{MR2925385}
Hanfeng Li.
\newblock Compact group automorphisms, addition formulas and
  {F}uglede-{K}adison determinants.
\newblock {\em Ann. of Math. (2)}, 176(1):303--347, 2012.

\bibitem[LL16]{li-liang-mean-length}
Hanfeng Li and Bingbing Liang.
\newblock Sofic mean length.
\newblock {\em preprint}, 2016.

\bibitem[LSW90]{LSW-1990}
Douglas Lind, Klaus Schmidt, and Tom Ward.
\newblock Mahler measure and entropy for commuting automorphisms of compact
  groups.
\newblock {\em Invent. Math.}, 101(3):593--629, 1990.

\bibitem[LT14]{MR3110799}
Hanfeng Li and Andreas Thom.
\newblock Entropy, determinants, and {$L^2$}-torsion.
\newblock {\em J. Amer. Math. Soc.}, 27(1):239--292, 2014.

\bibitem[Mal40]{malcev-1940}
A.~Malcev.
\newblock On isomorphic matrix representations of infinite groups.
\newblock {\em Rec. Math. [Mat. Sbornik] N.S.}, 8 (50):405--422, 1940.

\bibitem[Orn70a]{ornstein-1970a}
Donald Ornstein.
\newblock Bernoulli shifts with the same entropy are isomorphic.
\newblock {\em Advances in Math.}, 4:337--352 (1970), 1970.

\bibitem[Orn70b]{ornstein-1970c}
Donald Ornstein.
\newblock Factors of {B}ernoulli shifts are {B}ernoulli shifts.
\newblock {\em Advances in Math.}, 5:349--364 (1970), 1970.

\bibitem[OW80]{OW80}
Donald~S. Ornstein and Benjamin Weiss.
\newblock Ergodic theory of amenable group actions. {I}. {T}he {R}ohlin lemma.
\newblock {\em Bull. Amer. Math. Soc. (N.S.)}, 2(1):161--164, 1980.

\bibitem[OW87]{OW87}
Donald~S. Ornstein and Benjamin Weiss.
\newblock Entropy and isomorphism theorems for actions of amenable groups.
\newblock {\em J. Analyse Math.}, 48:1--141, 1987.

\bibitem[P{\u{a}}u11]{paunescu-2011}
Liviu P{\u{a}}unescu.
\newblock On sofic actions and equivalence relations.
\newblock {\em J. Funct. Anal.}, 261(9):2461--2485, 2011.

\bibitem[Pes08]{pestov-sofic-survey}
Vladimir~G. Pestov.
\newblock Hyperlinear and sofic groups: a brief guide.
\newblock {\em Bull. Symbolic Logic}, 14(4):449--480, 2008.

\bibitem[Pet89]{petersen-book}
Karl Petersen.
\newblock {\em Ergodic theory}, volume~2 of {\em Cambridge Studies in Advanced
  Mathematics}.
\newblock Cambridge University Press, Cambridge, 1989.
\newblock Corrected reprint of the 1983 original.

\bibitem[PK12]{pestov-kwiatkowska}
Vladimir~G Pestov and Alexsandra Kwiatkowska.
\newblock An introduction to hyperlinear and sofic groups.
\newblock {\em arXiv preprint arXiv:0911.4266}, 2012.

\bibitem[PS07]{popa-sasyk}
Sorin Popa and Roman Sasyk.
\newblock On the cohomology of {B}ernoulli actions.
\newblock {\em Ergodic Theory Dynam. Systems}, 27(1):241--251, 2007.

\bibitem[Roh67]{rohlin-lectures}
V.~A. Rohlin.
\newblock Lectures on the entropy theory of transformations with invariant
  measure.
\newblock {\em Uspehi Mat. Nauk}, 22(5 (137)):3--56, 1967.

\bibitem[Sch95]{schmidt-book}
Klaus Schmidt.
\newblock {\em Dynamical systems of algebraic origin}.
\newblock Modern Birkh\"auser Classics. Birkh\"auser/Springer Basel AG, Basel,
  1995.
\newblock [2011 reprint of the 1995 original] [MR1345152].

\bibitem[Sew14a]{seward-ergodic}
Brandon Seward.
\newblock Finite entropy actions of free groups, rigidity of stabilizers, and a
  howe--moore type phenomenon.
\newblock {\em To appear in Journal d'Analyse Math\'ematique}, 2014.

\bibitem[Sew14b]{seward-kreiger-1}
Brandon Seward.
\newblock Krieger's finite generator theorem for ergodic actions of countable
  groups {I}.
\newblock {\em arXiv:1405.3604}, 2014.

\bibitem[Sew14c]{seward-subgroup}
Brandon Seward.
\newblock A subgroup formula for f-invariant entropy.
\newblock {\em Ergodic Theory Dynam. Systems}, 34(1):263--298, 2014.

\bibitem[Sew15a]{seward-kreiger-2}
Brandon Seward.
\newblock Krieger's finite generator theorem for ergodic actions of countable
  groups {II}.
\newblock {\em arXiv:1501.03367v2}, 2015.

\bibitem[Sew15b]{seward-sinai-30}
Brandon Seward.
\newblock Positive entropy actions of countable groups factor onto {B}ernoulli
  shifts.
\newblock {\em preprint, version 30}, 2015.

\bibitem[STD16]{seward-tucker-drob}
Brandon Seward and Robin~D. Tucker-Drob.
\newblock Borel structurability on the 2-shift of a countable group.
\newblock {\em Ann. Pure Appl. Logic}, 167(1):1--21, 2016.

\bibitem[Ste75]{stepin-1975}
A.~M. Stepin.
\newblock Bernoulli shifts on groups.
\newblock {\em Dokl. Akad. Nauk SSSR}, 223(2):300--302, 1975.

\bibitem[Tho08]{MR2417890}
Andreas Thom.
\newblock Sofic groups and {D}iophantine approximation.
\newblock {\em Comm. Pure Appl. Math.}, 61(8):1155--1171, 2008.

\bibitem[Voi95]{MR1403924}
Dan Voiculescu.
\newblock Free probability theory: random matrices and von {N}eumann algebras.
\newblock In {\em Proceedings of the {I}nternational {C}ongress of
  {M}athematicians, {V}ol.\ 1, 2 ({Z}\"urich, 1994)}, pages 227--241.
  Birkh\"auser, Basel, 1995.

\bibitem[Wei00]{weiss-2000}
Benjamin Weiss.
\newblock Sofic groups and dynamical systems.
\newblock {\em Sankhy\=a Ser. A}, 62(3):350--359, 2000.
\newblock Ergodic theory and harmonic analysis (Mumbai, 1999).

\bibitem[Wei15]{MR3411529}
Benjamin Weiss.
\newblock Entropy and actions of sofic groups.
\newblock {\em Discrete Contin. Dyn. Syst. Ser. B}, 20(10):3375--3383, 2015.

\bibitem[Zha12]{zhang-local-variational}
Guohua Zhang.
\newblock Local variational principle concerning entropy of sofic group action.
\newblock {\em J. Funct. Anal.}, 262(4):1954--1985, 2012.

\end{thebibliography}
\bibliographystyle{alpha}

\end{document}